\documentstyle[amscd,amssymb,verbatim]{amsart}
\newcommand{\ad}{\operatorname{ad}}

\newcommand{\splin}{\operatorname{sl}}

\newcommand{\G}{{\Bbb G}}

\newcommand{\gog}{{\frak g}}

\newcommand{\Pic}{\operatorname{Pic}}

\newcommand{\Ad}{\operatorname{Ad}}

\numberwithin{equation}{section}

\newtheorem{thm}{Theorem}[section]
\newtheorem{prop}[thm]{Proposition}
\newtheorem{lem}[thm]{Lemma}
\newtheorem{cor}[thm]{Corollary}

\newenvironment{rem}{\vspace{3mm}\noindent
{\bf Remark.}}{\vspace{3mm}}
\newenvironment{rems}{\vspace{3mm}
\noindent {\bf Remarks.}}{\vspace{3mm}}

\newcommand{\Pf}{\noindent {\it Proof}}
\newcommand{\id}{\operatorname{id}}

\newcommand{\PP}{{\cal P}}

\newcommand{\End}{\operatorname{End}}

\newcommand{\QT}{{\cal QT}}

\newcommand{\De}{\Delta}

\newcommand{\C}{{\Bbb C}}

\newcommand{\Z}{{\Bbb Z}}
\newcommand{\Q}{{\Bbb Q}}

\newcommand{\La}{\Lambda}

\newcommand{\wt}{\widetilde}
\newcommand{\ot}{\otimes}

\newcommand{\sub}{\subset}
\newcommand{\ed}{\qed\vspace{3mm}}

\newcommand{\CH}{\operatorname{CH}}
\newcommand{\NS}{\operatorname{NS}}

\title{Fourier-stable subrings in the Chow rings of abelian varieties}
\author{A. Polishchuk}
\address{Department of Mathematics, University of Oregon, Eugene, OR 97403}
\email{apolish@@uoregon.edu}
\thanks{This work was partially supported by the NSF grant DMS-0601034}

\begin{document}
\begin{abstract} We study subrings in the Chow ring $\CH^*(A)_{\Q}$ of an abelian variety $A$,
stable under the Fourier transform with respect to an arbitrary polarization. We prove that
by taking Pontryagin products of classes of dimension $\le 1$ one gets such a subring.
We also show how to construct finite-dimensional Fourier-stable subrings in $\CH^*(A)_{\Q}$. Another
result concerns the relation between the Pontryagin product and the usual product on
the $\CH^*(A)_{\Q}$. We prove that the operator of the usual product with a cycle
is a differential operator with respect to the Pontryagin product and compute its order in terms
of the Beauville's decomposition of $\CH^*(A)_{\Q}$.
\end{abstract}
\maketitle

\bigskip

\centerline{\sc Introduction}

\bigskip

The goal of this paper is to generalize some facts about cycles on the Jacobian variety of a curve
observed in \cite{B} and \cite{P-lie} to the case of an arbitrary abelian variety $A$. 
Recall that the Fourier transform identifies the Chow group of $A$ with rational coefficients
$\CH^*(A)_{\Q}$ with $\CH^*(\hat{A})_{\Q}$, where $\hat{A}$ is the dual abelian variety
(see \cite{M}, \cite{B1}). If $A$ is equipped with a polarization $d$ then one can define
the Fourier transform $F_d$ which is an automorphism of $\CH^*(A)_{\Q}$.  
Our main result is concerned with constructing Fourier-stable subrings in $\CH^*(A)_{\Q}$.
Note that since the Fourier transform interchanges the usual intersection-product with the
Pontryagin product on $\CH^*(A)_{\Q}$, such subrings are closed under both products.
Beauville has shown in \cite{B} that starting with the class $[C]$
of a curve in its Jacobian $J$ and taking the subalgebra in $\CH^*(J)_{\Q}$ with respect to the Pontryagin product generated by the classes $[m]_*[C]$, $m\in\Z$, 
one obtains a Fourier-stable subring in $\CH^*(J)_{\Q}$,
finite-dimensional modulo algebraic equivalence (here $[m]$ denotes the endomorphism of multiplication by $m\in\Z$ on an abelian variety).
In general one can ask for which finite-dimensional subspaces 
$V\sub\CH^*(A)_{\Q}$ the subalgebra with respect to the
Pontryagin product generated by $V$ is Fourier-stable (with respect to one or all polarizations).
We will show that there are many such subspaces, and in particular,
the space of classes of dimension $\le 1$ is the union of subspaces with this property
(see Theorem \ref{finite-thm}).

One Fourier-stable subring in $\CH^*(A)_{\Q}$ seems to deserve a special attention
although it is in general infinite-dimensional (even modulo algebraic equivalence).
Namely, we define the {\it quasitautological subring} $\QT(A)\sub\CH^*(A)_{\Q}$ as
the subring with respect to the Pontryagin product generated by $\CH^{g}(A)_{\Q}$ and
$\CH^{g-1}(A)_{\Q}$, where $g=\dim A$. We show that these subrings are Fourier-stable,
contain all the divisorial classes (intersections
of divisors) and are stable under push-forwards and pull-backs with respect to homomorphisms
between abelian varieties (see Theorem \ref{taut-thm} and Proposition \ref{taut-prop}). 

Another result of this paper is inspired by the appearance of differential operators when considering
various natural operators on tautological cycles on the Jacobian (see \cite{P-lie}).
We prove that the intersection-product with any divisor class is a differential operator of order $\le 2$
on $\CH^*(A)$ with respect to the Pontryagin product. Furthermore, the symbol of this operator
depends only on the corresponding class in the N\'eron-Severi group.
The intersection-product with an arbitrary class is still a differential operator and we show how
to compute its order using the Beauville's decomposition of $\CH^*(A)_{\Q}$ (see Theorem \ref{diff-thm}).
These results hold also for abelian schemes.

Finally, we show that the symbol of the second order differential operator given by the product with
an ample divisor is related to the Jordan algebra structure on the N\'eron-Severi group of an abelian
variety (see Proposition \ref{Jordan-prop}).

\section{Differential operators coming from cycles}

Let $G$ be a commutative proper flat group scheme over $S$ (e.g., an abelian scheme).
Recall that the Pontryagin product on $\CH_*(G)$ is given by
$$x*y=m_*(p_1^*x\cdot p_2^*y),$$
where $m:G\times_S G\to G$ is the group law, $p_1,p_2:G\times_S G\to G$ are the projections.
Now let $\xi$ be a biextension of $G\times_S G$ by $\G_m$.
We think about $\xi$ as a line bundle over $G\times_S G$.
Then we can associate with $\xi$ the following symmetric binary operation on $\CH_*(G)$:
$$\{x,y\}_{\xi}=m_*(c_1(\xi)\cdot p_1^*x\cdot p_2^*y),$$
Recall that if $L$ is a line bundle on $G$ equipped with a cube structure (see \cite{Breen})
then $\xi=m^*L\ot p_1^*L^{-1}\ot p_2^*L^{-1}$ has a natural biextension structure.
We say in this case that $\xi$ is associated with the cube structure $L$.
This is especially useful if $G$ is an abelian scheme because in this case a cube structure is the 
same as a trivialization along the zero section.

\begin{prop}\label{biext-prop} 
For all $x,y,z\in\CH_*(G)$ one has 
\begin{equation}\label{bider-eq}
\{x*y,z\}_{\xi}=x*\{y,z\}_{\xi}+y*\{x,z\}_{\xi}.
\end{equation}
If $\xi$ is associated with a cube structure $L$ then
\begin{equation}\label{biext-id}
\{x,y\}_{\xi}=c_1(L)\cdot(x*y)-(c_1(L)\cdot x)*y-(c_1(L)\cdot y)*x.
\end{equation}
\end{prop}

\Pf . Using the projection formula we get
$$\{x*y,z\}_{\xi}=(p_1+p_2+p_3)_*((p_1+p_2,p_3)^*c_1(\xi)\cdot p_1^*x\cdot p_2^*y\cdot p_3^*z),$$
where $p_i:G\times_S G\times_S G\to G$ are the projections.
Since $\xi$ is a biextension, we have 
$$(p_1+p_2,p_3)^*c_1(\xi)=p_{13}^*c_1(\xi)+p_{23}^*c_1(\xi).$$
This immediately leads to formula \eqref{bider-eq}.
Substituting the equality $c_1(\xi)=m^*c_1(L)-p_1^*c_1(L)-p_2^*c_1(L)$
into the definition of
$\{x,y\}_{\xi}$ we derive \eqref{biext-id}.
\ed

\begin{cor}\label{cube-cor}
For every line bundle with a cube structure $L$ on $G$ the operator
of multiplication by $c_1(L)$ on $\CH_*(G)$ is a differential operator
of order $\le 2$ with respect to the Pontryagin product.
\end{cor}

Now we are going to specialize to the case of abelian schemes. We assume that the base
scheme $S$ is smooth quasiprojective over a field.
Recall that for an abelian scheme $A/S$ there is a canonical decomposition
$$\CH^p(A)_{\Q}=\oplus_{s}\CH^p_s(A),$$
where $[m]^*x=m^{2p-s}x$ for $x\in\CH^p_s(A)$, $m\in\Z$ (see \cite{B2}, \cite{DM}).
Moreover, the component $\CH^p_s(A)$ is nonzero only for $0\le 2p-s\le 2g$, where
$g$ is the relative dimension of $A/S$.

By a {\it polarization} of an abelian scheme $A/S$ we mean a relatively ample
class $d\in\Pic(A)_{\Q}$ such that $[0]^*d=0$ (i.e., $d$ is trivialized along the zero section) and 
$[-1]^*d=d$ ($d$ is symmetric). Note that these two last conditions are equivalent to $d\in\CH^1_0(A)$.
Recall that with such a polarization one can associate a
Lefschetz $\splin_2$-action on $\CH^*(A)_{\Q}$ as follows (see \cite{K}): 
$$e(x)=d\cdot x,\ \ f(x)=\frac{d^{g-1}}{(g-1)!\chi(d)}*x, \ \ h|_{\CH^p_s(A)}=(2p-s-g)\id,$$
where $\chi(d)$ is the square root of the degree of the isogeny $\phi:A\to\hat{A}$ associated
with $d$.

Let $F$ denote the Fourier transform $F:\CH^*(A)_{\Q}\to\CH^*(\hat{A})_{\Q}$ (see \cite{B1}).
For a polarization $d$ we set 
$$F_d=\frac{1}{\chi(d)}\phi^*\circ F:\CH^*(A)_{\Q}\to\CH^*(A)_{\Q}.$$
Then one has $F_d^2=(-1)^g[-1]^*$ (see \cite{K}, Lemma 6.1).
Also, $F_d$ exchanges the Pontryagin product with the usual product up to a constant:
\begin{equation}\label{four-prod-eq}
F_d(x*y)=\chi(d)\cdot F_d(x)\cdot F_d(y).
\end{equation}
The Fourier transform intertwines the above $\splin_2$-action in the following way:
\begin{equation}\label{inter-eq}
F_deF_d^{-1}=-f, \ \ F_dfF_d^{-1}=-e, \ \ F_dhF_d^{-1}=-h.
\end{equation}

\begin{cor}\label{sl2-diff-cor} 
Let $A/S$ be a polarized abelian scheme.
Then the operator $e$ (resp., $f$) of the associated Lefschetz $\splin_2$-action on $\CH^*(A)_{\Q}$
is a differential operator of order $\le 2$
with respect to the Pontryagin (resp., usual) product.
\end{cor}

\Pf . The assertion about $e$ follows from Corollary \ref{cube-cor} since $d$ is trivialized along the zero
section. The assertion about $f$ follows by Fourier duality.
\ed

Using the $\splin_2$-action associated with a polarization we will prove a much more general
statement in Theorem \ref{diff-thm} below. The following lemma in the principally polarized case
is essentially equivalent to the identity (1.7) of \cite{B}, and our proof is an easy adaptation
of the same argument. 

\begin{lem}\label{exp-lem} 
The operators $F_d$, $e$ and $f$ associated with a polarization $d$ on an abelian scheme
$A/S$ satisfy
$$(-1)^gF_d=\exp(e)\exp(-f)\exp(e).$$
\end{lem}

\Pf . Let $p_1,p_2$ denote the projections from the product $A\times_S\hat{A}$ to the factors.
Then
$$\chi(d)F_d(x)=\phi^*p_{2*}(e^{c_1(\PP)}\cdot p_1^*x),$$
where $\PP$ is the Poincar\'e line bundle on $A\times_S\hat{A}$.
Since $(\id\times\phi)^*\PP$ is the biextension associated with $d$, we have
$$(\id\times\phi)^*c_1(\PP)=m^*d-p_1^*d-p_2^*d.$$
Hence,
$$\chi(d)F_d(x)=e^{-d}\cdot p_{2*}(p_1^*(e^{-d}\cdot x)\cdot e^{m^*d}).$$
Making the change of variables $(x,y)\mapsto (-x,x+y)$ on $A\times A$, we can rewrite this as
\begin{equation}\label{exp-four-aux-eq}
\chi(d)F_d(x)=e^{-d}\cdot m_*(p_1^*(e^{-d}\cdot [-1]^*x)\cdot p_2^*e^d)=
e^{-d}\cdot[(e^{-d}\cdot [-1]^*x)*e^d].
\end{equation}
Since $F_d(e^d)=e^{-d}$ (see \cite{K}, Prop. 2.2), applying \eqref{four-prod-eq} we get
$$e^d*y=\chi(d) F_d^{-1}(F_d(e^d)\cdot F_d(y))=\chi(d)F_d^{-1}(e^{-d}\cdot F_d(y))=
\chi(d)F_d^{-1}\exp(-e)F_d(y)=\chi(d)\exp(f)y.$$
Using this for $y=e^{-d}\cdot [-1]^*x$ we can rewrite \eqref{exp-four-aux-eq} as
$$F_d[-1]^*=\exp(-e)\exp(f)\exp(-e).$$
Passing to inverses we get the required identity.
\ed

\begin{thm}\label{diff-thm}
Let $A/S$ be an abelian scheme of relative dimension $g$. 
For every nonzero class $a\in\CH^p_s(A)$ the operator $x\mapsto a\cdot x$ (resp.,
$x\mapsto a*x$) is a differential operator on $\CH^*(A)_{\Q}$
with respect to the Pontryagin product
(resp., usual product) of order $2p-s$ (resp., $2g-2p+s$).
\end{thm}

For a class $a\in\CH^*(A)_{\Q}$ let us denote the operators of the usual (resp., Pontryagin) product
with $a$ as follows:
$$L_a(x)=a\cdot x, \ \ \La_a(x)=a*x.$$
For an operator $B$ on a vector space $V$ we denote by $\ad(B)$ the 
corresponding operator $X\mapsto BX-XB$ on the algebra of endomorphisms of $V$.

\begin{lem}\label{sl2-lem}
Let $d\in\CH^1_0(A)$ be a polarization.
Then for $a\in\CH^p_s(A)$ one has
$$L_{F_d(a)}=c\cdot \ad(e)^{2g-2p+s}(\La_a),$$
where $c\in\Q^*$.
\end{lem}

\Pf . Let us consider the adjoint $\splin_2$-action on $\End(\CH^*(A))$. Then we have
$\ad(f)(\La_a)=0$, $\ad(h)(\La_a)=(2p-s-2g)\La_a$. Therefore, $\La_a$ is the lowest
weight vector of weight $-(2g-2p+s)$ (recall that $2g-2p+s\ge 0$). 
Note that from \eqref{four-prod-eq} we get the equality of operators
$$L_{F_d(a)}=\chi(d)^{-1}\cdot \Ad(F_d)(\La_a),$$
where for $X\in\End(\CH^*(A))$ we set $\Ad(F_d)(X)=F_dXF_d^{-1}$.
By Lemma \ref{exp-lem},
the operator $\Ad(F_d)$ preserves the $\splin_2$-submodule in $\End(\CH^*(A))$
generated by $\La_a$. Furthermore, since $\Ad(F_d)$ intertwines the adjoint action of $\splin_2$
in a way similar to \eqref{inter-eq}, it exchanges the lowest and highest weight components
in this irreducible representation. Since $\ad(e)^{2g-2p+s}(\La_a)$ generates the highest weight
component, this implies the result.
\ed

\noindent
{\it Proof of Theorem \ref{diff-thm}.}
Let us fix a polarization $d$ on $A/S$.
By Corollary \ref{sl2-diff-cor}, the corresponding operator $e$ is a
differential operator of order $\le 2$ with respect to the Pontryagin product.
By Lemma \ref{sl2-lem}, this implies that for $a\in\CH^p_s(A)$ 
the operator $L_{F_d(a)}$ is of order $\le 2g-2p+s$ (taking the 
commutator with a second order differential operator raises the order by one).
Since $F_d$ exchanges $\CH^p_s(A)$ and $\CH^{g-p+s}_s(A)$,
this shows that for every $b\in\CH^q_s(A)$
the operator $L_b$ has order $\le 2q-s$ with respect to the Pontryagin product.
By Fourier duality, it follows that for $a\in\CH^p_s(A)$ the operator $\La_a$ has order
$\le 2g-2p+s$ with respect to the usual product. Using the identity of Lemma \ref{sl2-lem}
again, we see that for $a\neq 0$ this order is exactly $2g-2p+s$ (recall that the operator $e$ is given by
the usual product with the class $d$).
\ed

\begin{rem} Using Theorem \ref{diff-thm} we can restate the Beauville's conjecture on the
vanishing of $\CH^p_s(A)$ for $s<0$ (see \cite{B2}) as follows:

\noindent {\it For every $a\in\CH^p(A)$ the differential operator $L_a$ (with respect to the Pontryagin
product) has order $\le 2p$.}
\end{rem}

\section{Fourier-stable subrings}

\begin{thm}\label{taut-thm} 
Let $A$ be an abelian variety over a field $k$.
Let us denote by $\QT^*(A)\sub\CH^*(A)_{\Q}$ the subring with respect to
the Pontryagin product generated by $\CH^{g}(A)_{\Q}$ and $\CH^{g-1}(A)_{\Q}$.
Then $\QT^*(A)$ is stable under the usual product and under
the Fourier transform with respect to any polarization of $A$. 
\end{thm}

This is an easy consequence of the following result.

\begin{lem}\label{Pon-lem} 
Let $d\in\CH^1_0(A)$ be a polarization, and let $\xi$
be the corresponding biextension on $A\times A$. 
Suppose $V\sub\CH^{g}(A)_{\Q}\oplus\CH^{g-1}(A)_{\Q}$
is a subspace closed under $\{\cdot,\cdot\}_{\xi}$, 
such that $d\cdot V\sub V$ and $[m]^*V\sub V$ for all $m\in\Z$. 
Then the subalgebra $R\sub\CH^*(A)_{\Q}$ with respect to the Pontryagin product 
generated by $\wt{V}=V+\Q d^{g-1}$ over $\Q[0]$
is invariant under the Fourier transform $F_d$ and under the usual product. 
\end{lem}

\Pf . Let us consider the $\splin_2$-action on $\CH^*(A)_{\Q}$ associated with $d$.
By assumption, we have $e(V)\sub V$ and hence $e(\wt{V})\sub \Q[0]+V\sub R$. 
Since $V$ respects Beauville's decomposition, we also have $h(V)\sub V$. Recall that
$f$ acts by the Pontryagin product with $\frac{d^{g-1}}{(g-1)!\chi(d)}$. Therefore,
$ef(\wt{V})\sub R$ (using the identity $ef-fe=h$).
Next, using \eqref{biext-id} and our assumption that $V$ is closed under $\{\cdot,\cdot\}_{\xi}$, 
we obtain $e(V*V)=d\cdot (V*V)\sub R$. 
Together with the inclusion $e(d^{g-1}*\wt{V})\sub R$ this implies that $e(\wt{V}*\wt{V})\sub R$.
Since $e$ is a differential operator of order $2$, we derive that $e(R)\sub R$.
Also, $f(R)\sub R$ by the definition of $R$. Therefore, $R$ is preserved by the $\splin_2$-action.
Hence, it is also stable under the Fourier transform 
$F_d$ (by Lemma \ref{exp-lem}). Now \eqref{four-prod-eq}
implies that $R$ is also closed under the usual product.
\ed

Note that in the situation of Theorem \ref{taut-thm} we take $V$ to be $\CH^{g}(A)_{\Q}\oplus
\CH^{g-1}(A)_{\Q}$, so the assumptions of the above lemma are satisfied for trivial reasons.

We will call $\QT^*(A)\sub\CH^*(A)_{\Q}$ the {\it quasitautological subring}.
Note that by definition, $\QT^*(A)$ is the $\Q$-linear span of $0$-cycles and of
classes of the form
$(f_1,\ldots,f_n)_*[C_1\times\ldots\times C_n]$,
where $f_i:C_i\to A$ are morphisms from curves to $A$ and
$(f_1,\ldots,f_n)(x_1,\ldots,x_n)=f(x_1)+\ldots+f(x_n)$. 
Let us list some other properties of this subring.

\begin{prop} \label{taut-prop} 
(i) $\QT^*(A)_{\Q}$ contains 
the subring generated by all the divisor classes (with respect to the usual product).

\noindent (ii) 
If $f:A\to B$ is a homomorphism of abelian varieties then 
$f_*\QT^*(A)_{\Q}\sub\QT^*(B)$ and $f^*\QT^*(B)\sub\QT^*(A)$.

\noindent (iii) $\QT^*(A)$ is a graded subspace with respect to the Beauville's decomposition. 
Also, $\QT^*(A)\sub\oplus_{s\ge 0}\CH^*_s(A)$.

\noindent (iv) The intersection $\QT^*(A)\cap\oplus_p\CH^p_0(A)$ consists of divisorial classes.
It $k=\C$ then the image of $\QT^*(A)$ in the cohomology algebra $H^*(A,\Q)$ coincides with
the subalgebra generated by the algebraic part of $H^2(A,\Q)$.
\end{prop}

\Pf . (i) We have seen in the proof of Theorem \ref{taut-thm} that $\QT^*(A)$ is stable under
the product with any divisor class. Also, since it is stable under the Fourier transform with respect
to some polarization $d$, it contains the class $[A]=F_d([0])$. Therefore, it contains all products
of divisors.

\noindent (ii) The assertion about $f_*$ is clear, since it is a homomorphism with respect to the Pontryagin product. The second assertion follows by Fourier duality.

\noindent (iii) The first assertion 
follows from the fact that $\QT^*(A)$ is stable under all operations $[m]_*$,
where $m\in\Z$. The second is implied by the vanishing $\CH^g_s(A)=\CH^{g-1}_s(A)=0$ for
$s<0$ (see \cite{B2}, Prop. 3).

\noindent (iv) The intersection in question is a Fourier-stable subalgebra.
Now to prove the first assertion we use the fact that 
the Fourier transform exchanges $\CH^{g-1}_0(A)$ and $\CH^1_0(A)$. The second assertion
follows from the first.
\ed

For example, for $g=3$ we have $\QT^*(A)=\CH^*(A)_{\Q}$. For $g=4$ the subring $\QT^*(A)$
is almost the entire $\CH^*(A)_{\Q}$. Namely, this is true for all the summands of Beauville's decomposition except for $\CH^2_0(A)$ because the intersection $\QT^*(A)\cap\CH^2_0(A)$
consists only of divisorial classes.

Note that in the case of the Jacobian the quasitautological subring is in general larger than the tautological subring defined by Beauville. This happens already for generic abelian threefold since in
this case $\CH^{2}(A)_{\Q}$ is infinite-dimensional modulo algebraic equivalence (see \cite{Nori}).

Slightly generalizing the idea of Lemma \ref{Pon-lem} 
we can construct a large class of Fourier-stable finite-dimensional
subrings of $\CH^*(A)_{\Q}$. 

\begin{thm}\label{finite-thm} 
Let $A$ be an abelian variety of dimension $g$ over a field $k$. 
Then for every finite-dimensional subspace $V\sub\CH^*_{>0}(A)=\oplus_{p\ge 0,s>0}\CH^p_s(A)$ 
there exists a finite-dimensional bigraded subspace $\wt{V}\sub\CH^*_{>0}(A)$ containing $V$
such that for $W=\wt{V}\oplus \CH^{g-1}_0(A)$ the subalgebra (with respect to the Pontryagin product)
$$R=\Q[0]+W+W*W+...$$ 
satisfies the following properties:

\noindent
(a)  $F_d(R)\sub R$ for any polarization $d$;

\noindent
(b) $R$ is a subring with respect to the usual product.

If we start with $V\sub\CH^{g}(A)_{\Q}\oplus \CH^{g-1}(A)_{\Q}$ 
then we can choose $\wt{V}\sub\CH^{g}(A)_{\Q}\oplus \CH^{g-1}(A)_{\Q}$ with the above properties.
\end{thm}

\Pf . Without loss of generality we can assume that $V$ is a bigraded subspace.
Let us denote by $\gog\sub\End_{\Q}(\CH^*(A)_{\Q})$ the Lie subalgebra generated by
all the $\splin_2$-triples associated with polarizations on $A$. Note that all these
$\splin_2$-triples have the common operator $h$. The adjoint action of $h$ gives $\gog$ a natural
grading. In fact, it is not difficult to see that the action of $\gog$ integrates to an algebraic representation 
of a reductive algebraic group on $\CH^*(A)_{\Q}$, so that every vector is contained in a finite-dimensional subrepresentation (see \cite{P-thesis}, Theorem 13.1).

Replacing $V$ with $U(\gog_{\ge 0})V$
we can assume that $U(\gog_{\ge 0})V\sub V$.
Next, we claim that we can embed $V$ into
a finite-dimensional subspace $\wt{V}\sub\CH^*_{>0}(A)$ such that
$U(\gog_{\ge 0})\wt{V}\sub \wt{V}$ and $U(\gog_{>0})_+(\wt{V}*\wt{V})\sub\wt{V}$,
where $U(\gog_{>0})_+\sub U(\gog_{>0})$ is the augmentation ideal.

We want to keep track only of the second grading on $\CH^*(A)_{\Q}$, so we denote
$V=\oplus_{s>0}V_s$, where $V_s\sub\oplus_p\CH^p_s(A)$. Note that all operators
in $\gog$ preserve this grading.
We will construct $\wt{V}$ by iterating the following procedure.
At each iteration we start with a subspace $V=\oplus_{s>0}V_s\sub\CH^*_{>0}(A)$ 
closed under the action of $U(\gog_{\ge 0})$ and replace it with a bigger $U(\gog_{\ge 0})$-submodule
$V'=V+U(\gog_0)U(\gog_{> 0})_+(V*V)$.
Let us show that the obtained sequence of subspaces stabilizes after a finite number of steps.
Indeed, assume that after some number of steps the input subspace satisfies
$U(\gog_{>0})_+(V_s*V_t)\sub V$ 
for $s+t\le n$ (initially this condition is satisfied for $n=1$).
We claim that after the next step the same condition will hold for $s+t\le n+1$. Indeed,
we have 
$$V'_r=V_r+\sum_{s+t=r}U(\gog_0)U(\gog_{>0})_+(V_s*V_t).$$ 
The above assumption implies that for $r\le n$ we have 
$U(\gog_0)U(\gog_{>0})_+(V_s*V_t)\sub V$ whenever
$s+t=r$ (recall that $V$ is closed under the action of $U(\gog_{\ge 0})$).
Hence, $V'_r=V_r$ for $r\le n$.
Therefore, for positive $s$ and $t$ such that $s+t\le n+1$ we have
$$U(\gog_{>0})_+(V'_s*V'_t)=U(\gog_{>0})_+(V_s*V_t)\sub V',$$
which proves our claim.
Thus, after a finite number of steps we will get a subspace $\wt{V}$ with the required
properties. Furthermore, if we started with a subspace in $\CH^{g}(A)_{\Q}\oplus\CH^{g-1}(A)_{\Q}$
then we will still have $\wt{V}\sub\CH^{g}(A)_{\Q}\oplus\CH^{g-1}(A)_{\Q}$.

We claim that the subalgebra $R\sub\CH^*(A)_{\Q}$ with respect to the Pontryagin product
generated by $W=\wt{V}\oplus\CH^{g-1}_0(A)$
is closed under the product with any class in $\CH^1_0(A)$.
Let us choose a basis $d_1,\ldots,d_N$ in $\CH^1_0(A)$ consisting of
ample classes, so that for every $i=1,\ldots,N$ we have the corresponding
$\splin_2$-action on $\CH^*(A)_{\Q}$
with $e_i(x)=d_i\cdot x$ and $f_i(x)=c_i*x$.
Since each $e_i$ is a differential operator of order $\le 2$ with respect to the Pontryagin product,
it is enough to check that $e_i(W)\sub R$ and $e_i(W*W)\sub R$ for every $i$.
Recall that by the construction we have $e_i(\wt{V})\sub\wt{V}$ and
$e_i(\wt{V}*\wt{V})\sub\wt{V}$.
Since $e_i(\CH^{g-1}_0(A))\sub\Q[0]$ and $e_i(\CH^{g-1}_0(A)*\CH^{g-1}_0(A))\sub\CH^{g-1}_0(A)$,
it remains to check that $e_i(\CH^{g-1}_0(A)*\wt{V})\sub R$.
The classes $(c_i)_{1\le i\le N}$ span $\CH^{g-1}_0(A)$, so the assertion is equivalent to
$e_if_j(\wt{V})\sub R$. But this follows from the fact that $[e_i,f_j](\wt{V})\sub\wt{V}$ (since
$[e_i,f_j]\in\gog_0$) and from the inclusion $f_je_i(\wt{V})\sub f_j(\wt{V})=c_j*\wt{V}\sub R$.

Since $\CH^1_0(A)\cdot R\sub R$ and $\CH^{g-1}_0(A)*R\sub R$, we see that $R$ is stable
under the $\splin_2$-action associated with any polarization of $A$. Hence, it is also
stable under the corresponding Fourier transforms. It follows that $R$ is stable under the usual
product.
\ed

We end with several observations about the
operation $\{\cdot,\cdot\}_{\xi}$ on $\CH^{g-1}(A)$ that played a role in Lemma \ref{Pon-lem}.
Let us introduce some notation. We set $\End^0(A)=\End(A)\ot\Q$, $\NS^0(A)=\CH^1_0(A)$. 
It is well-known that a choice of polarization $d$ on $A$ gives rise to an isomorphism
$$\rho:\NS^0(A)\wt{\to} \End^0(A)^+: [L]\mapsto \phi_d^{-1}\circ\phi_L,$$
where $\phi_L:A\to\hat{A}$ denotes the symmetric homomorphism associated with a line bundle $L$,
$\End^0(A)^+\sub\End^0(A)$ is the subspace of elements
invariant under the Rosati involution associated with $d$ (sending $f$ to $\phi_d^{-1}\hat{f}\phi_d$).
We denote by $f\mapsto L(f)$ the inverse map to $\rho$. In fact, one has 
$L(f)=\frac{1}{2}(\id\times(\phi_d\circ f))^*c_1(\PP)$, where $\PP$ is Poincar\'e line bundle on 
$A\times\hat{A}$.
Thus, the map $f\mapsto F_d(L(f))$ is an isomorphism of $\End^0(A)^+$ onto $\CH^{g-1}_0(A)$.
Part (i) of the following proposition shows that under this isomorphism the operation
$\{\cdot,\cdot\}_{\xi}$ becomes the usual Jordan multiplication on $\End^0(A)^+$ (up to a constant).

\begin{prop}\label{Jordan-prop} 
(i) Let $\xi$ be the biextension associated with a polarization $d\in\CH^1_0(A)_{\Q}$.
Then for $f_1,f_2\in\End^0(A)^+$ one has
$$\{F_d(L(f_1)),F_d(L(f_2))\}_{\xi}=(-1)^g\chi(d) F_d(L(f_1f_2+f_2f_1)).$$

\noindent (ii) In the above situation assume that
$x\in\CH^{g-1}_0(A)$ corresponds to an endomorphism of $A$ satisfying a quadratic equation
over $\Q$. Then for every $y\in\CH^{g-1}(A)_{\Q}$ 
the Jordan identity is satisfied:
$$\{\{x,y\}_{\xi},\{x,x\}_{\xi}\}_{\xi}=\{x,\{y,\{x,x\}_{\xi}\}_{\xi}\}_{\xi}.$$

\noindent (iii) Let $J$ be the Jacobian of a smooth projective curve $C$, $\xi$ be the biextension
corresponding to the standard principal polarization of $J$. 
We consider the embedding $C\to J$ sending a point
$p_0\in C$ to $0$, so that $\xi|_{C\times C}=-[\De]+[p_0\times C]+[C\times p_0]$. Then
$$\{[C]_s,[C]_t\}_{\xi}=-{s+t+2\choose s+1}[C]_{s+t},$$
where $[C]=\sum_{s\ge 0}[C]_s$ with $[C]_s\in\CH^{g-1}_s(J)$.
\end{prop}

\Pf . (i) We will prove a more general formula
\begin{equation}\label{J-gen-for}
\{F_d(e^{L(f_1)}),F_d(e^{L(f_2)})\}_{\xi}=(-1)^g\chi(d) F_d(L(f_1f_2+f_2f_1)\cdot e^{L(f_1+f_2)})
\end{equation}
from which the required identity is obtained by considering parts of codimension $g-1$.
Let us set for $f\in\End^0(A)^+$
$$N(f)=\chi(L(f))/\chi(d),$$
where $\chi:\NS^0(A)\to\Q$ is the polynomial function given by the Euler characteristic. 
Note that $N(f)^2=\deg(f)$, and hence, $N(f_1f_2)=N(f_1)N(f_2)$. 
Also, by Serre duality we have $N(-f)=(-1)^g N(f)$.
We will work with elements of a Zariski open subset of $\End^0(A)^+$, so that the inverses
of various elements are well-defined. Combining
the well-known formula for the Fourier transform of $e^{L(f)}$ (see e.g., \cite{B1}) with the identity
$(f^{-1})^*L(f)=L(f^{-1})$ for $f\in\End^0(A)^+$ we obtain
\begin{equation}\label{four-exp-eq}
F_d(e^{L(f)})=N(f)e^{L(-f^{-1})}.
\end{equation}
The identity \eqref{biext-id} implies that $\{x,y\}_{\xi}$ is the coefficient of $t$ in
$$e^{td}\cdot ((e^{-td}\cdot x)*(e^{-td}\cdot y)).$$
Therefore, we can deduce \eqref{J-gen-for} from the following equation:
\begin{equation}\label{four-grand-eq}
F_d^{-1}\left(e^{td}\cdot [(e^{-td}\cdot F_d(e^{L(f_1)}))*(e^{-td}\cdot F_d(e^{L(f_2)}))]\right)=
(-1)^g\chi(d)N(1+tf_1)N(1+tf_2)N(1-tf)e^{L(\frac{f}{1-tf})},
\end{equation}
where $f=\frac{f_1}{1+tf_1}+\frac{f_2}{1+tf_2}$. Indeed, the coefficient of $t$ can be easily extracted
since 
$$f=f_1+f_2-t(f_1^2+f_2^2)+\ldots,\ \ \frac{f}{1-tf}=f_1+f_2+t(f_1f_2+f_2f_1)+\ldots.$$
It remains to prove \eqref{four-grand-eq}. This is a straightforward calculation.
Note that $L(1)=d$. Therefore, using \eqref{four-exp-eq} and \eqref{four-prod-eq} we can write
\begin{align*}
& F_d\left((e^{-td}\cdot F_d(e^{L(f_1)}))*(e^{-td}\cdot F_d(e^{L(f_2)}))\right)=
N(f_1)N(f_2)F_d(e^{L(-t-f_1^{-1})}*e^{L(-t-f_2^{-1})})=\\
&N(f_1)N(f_2)\chi(d)F_d(e^{L(-t-f_1^{-1})})\cdot F_d(e^{L(-t-f_2^{-1})})=
\chi(d)N(1+tf_1)N(1+tf_2)e^{L(f)}.
\end{align*}
Applying $F_d$ again we deduce
$$(e^{-td}\cdot F_d(e^{L(f_1)}))*(e^{-td}\cdot F_d(e^{L(f_2)}))=(-1)^g\chi(d)N(1+tf_1)N(1+tf_2)N(f)
e^{L(-f^{-1})}.$$
Hence, the left-hand side of \eqref{four-grand-eq} is equal to
$$(-1)^g\chi(d)N(1+tf_1)N(1+tf_2)N(f)F_d^{-1}(e^{L(t-f^{-1})})=
(-1)^g\chi(d)N(1+tf_1)N(1+tf_2)N(1-tf)e^{L(\frac{f}{1-tf})},$$
as required.

\noindent (ii) By part (i) in this case $\{x,x\}_{\xi}\in\CH^{g-1}_0(A)$ 
is a linear combination of $x$ and $d^{g-1}$.
Since $\{d^{g-1},\cdot\}_{\xi}$ is proportional to $h-g\id$ and the operator $\{x,\cdot\}_{\xi}$
preserves the grading given by $h$, it follows that the operators 
$\{\{x,x\}_{\xi},\cdot\}_{\xi}$ and $\{x,\cdot\}_{\xi}$ commute.

\noindent (iii) Using the formula for the restriction of $\xi$ to $C\times C$ together with the isomorphism
$([m]\times[n])^*\xi=mn\xi$ we obtain
$$\{[m]_*[C],[n]_*[C]\}_{\xi}=-mn([m+n]_*[C]-[m]_*C-[n]_*C).$$
Taking into account the formula $[m]_*[C]=\sum_{s\ge }m^{s+2}[C]_s$ we get the result.
\ed

\begin{rems} 1. The first two parts of the above proposition (and their proofs) 
work for an abelian scheme as well. In the case of an abelian variety over $\C$ one can also
work in the cohomology ring and use the action of the corresponding N\'eron-Severi
algebra (see \cite{LL}).

\noindent 2.
Part (iii) of the above proposition shows that on the entire $\CH^{g-1}(A)_{\Q}$ the operation
$\{\cdot,\cdot\}_{\xi}$ in general does not satisfies the Jordan identity (and that
$\CH^{g-1}(A)_{\Q}$ is not a Jordan module over $\CH^{g-1}_0(A)$).
\end{rems}

\end{document}